\documentclass[a4paper,12pt]{amsart}
\usepackage{amssymb}
\usepackage{ifthen}
\usepackage{graphicx}
\usepackage{float}
\usepackage{caption}
\usepackage{subcaption}
\usepackage{cite}
\usepackage{amsfonts}
\usepackage{amscd}
\usepackage{amsxtra}
\usepackage{color}
\newcommand{\floor}[1]{\lfloor #1 \rfloor}
\newtheorem{thm}{Theorem}
\newtheorem{cor}{Corollary}
\newtheorem{lem}{Lemma}
\newtheorem{rem}{Remark}

\newtheorem{conj}{Conjecture}
\newtheorem{prob}{Problem}

\theoremstyle{definition}
\newtheorem{defn}{Definition}[section]
\newtheorem{example}{Example}

\newenvironment{pf}[1][]{%
 \vskip 1mm
 \noindent
 \ifthenelse{\equal{#1}{}}%
  {{\slshape Proof. }}%
  {{\slshape #1.} }%
 }%
{\qed\bigskip}

\newcounter{alphabet}

\newenvironment{Thm}[1][]{\refstepcounter{alphabet}%
\bigskip%
\noindent%
{\bf Theorem \Alph{alphabet}}%
\ifthenelse{\equal{#1}{}}{}{ (#1)}%
{\bf .} \itshape}{\vskip 8pt}
\newenvironment{Lem}[1][]{\refstepcounter{alphabet}%
\bigskip%
\noindent%
{\bf Lemma \Alph{alphabet}}%
{\bf .} \itshape}{\vskip 8pt}

\newcommand{\IN}{{\mathbb N}}
\newcommand{\IC}{{\mathbb C}}
\newcommand{\ID}{{\mathbb D}}

\def\be{\begin{equation}}
\def\ee{\end{equation}}

\newcommand{\bee}{\begin{enumerate}}
\newcommand{\eee}{\end{enumerate}}

\newcommand{\blem}{\begin{lem}}
\newcommand{\elem}{\end{lem}}
\newcommand{\bthm}{\begin{thm}}
\newcommand{\ethm}{\end{thm}}
\newcommand{\bcor}{\begin{cor}}
\newcommand{\ecor}{\end{cor}}
\newcommand{\beg}{\begin{example}}
\newcommand{\eeg}{\end{example}}
\newcommand{\begs}{\begin{examples}}
\newcommand{\eegs}{\end{examples}}
\newcommand{\bdefe}{\begin{defn}}
\newcommand{\edefe}{\end{defn}}
\newcommand{\bprob}{\begin{prob}}
\newcommand{\eprob}{\end{prob}}
\newcommand{\bques}{\begin{ques}}
\newcommand{\eques}{\end{ques}}
\newcommand{\bei}{\begin{itemize}}
\newcommand{\eei}{\end{itemize}}
\newcommand{\bcon}{\begin{conj}}
\newcommand{\econ}{\end{conj}}
\newcommand{\bcons}{\begin{conjs}}
\newcommand{\econs}{\end{conjs}}
\newcommand{\bprop}{\begin{propo}}
\newcommand{\eprop}{\end{propo}}
\newcommand{\br}{\begin{rem}}
\newcommand{\er}{\end{rem}}
\newcommand{\brs}{\begin{rems}}
\newcommand{\ers}{\end{rems}}
\newcommand{\bo}{\begin{obser}}
\newcommand{\eo}{\end{obser}}
\newcommand{\bos}{\begin{obsers}}
\newcommand{\eos}{\end{obsers}}
\newcommand{\bpf}{\begin{pf}}
\newcommand{\epf}{\end{pf}}
\newcommand{\ba}{\begin{array}}
\newcommand{\ea}{\end{array}}
\newcommand{\beq}{\begin{eqnarray}}
\newcommand{\beqq}{\begin{eqnarray*}}
\newcommand{\eeq}{\end{eqnarray}}
\newcommand{\eeqq}{\end{eqnarray*}}

\newcommand{\ds}{\displaystyle}

\newcounter{minutes}
\divide\time by 60
\newcounter{hours}
\multiply\time by 60 \addtocounter{minutes}{-\time}
\begin{document}

\bibliographystyle{amsplain}

%

\title[Modifications of Bohr's inequality in various settings]{Modifications of Bohr's inequality in various settings}

\thanks{
File:~\jobname .tex,
          printed: \number\day-\number\month-\number\year,
          \thehours.\ifnum\theminutes<10{0}\fi\theminutes}

\author[S. Ponnusamy]{Saminathan Ponnusamy
}
\address{
S. Ponnusamy, Department of Mathematics,
Indian Institute of Technology Madras, Chennai-600 036, India.
}
\email{samy@iitm.ac.in}

\author[R. Vijayakumar]{Ramakrishnan Vijayakumar}
\address{
R. Vijayakumar, Department of Mathematics,
Indian Institute of Technology Madras, Chennai-600 036, India.
}
\email{mathesvijay8@gmail.com}

\author[K.-J. Wirths]{Karl-Joachim Wirths}
\address{K.-J. Wirths, Institut f\"ur Analysis und Algebra, TU Braunschweig,
38106 Braunschweig, Germany.}
\email{kjwirths@tu-bs.de}

\subjclass[2010]{Primary: 30A10, 30B10 31A05, 30H05,  41A58; Secondary:  30C62, 30C75, 40A30
}
\keywords{Analytic functions, harmonic function, quasiconformal mapping, Bohr's inequality, subordination and quasisubordination
}

\begin{abstract}
The concept of Bohr radius for the class of bounded analytic functions was introduced by Harald Bohr in 1914.
His initial result received greater interest and was sharpened-refined-generalized by several mathematicians
in various settings--which is now called Bohr phenomenon. Various generalization of Bohr's classical theorem is now
an active area of research and has been a source of investigation in numerous
other function spaces and including holomorphic functions of several complex variables. Recently,
a new generalization of Bohr's ideas was introduced and investigated by Kayumov et al.. In this note,
we investigate and refine generalized Bohr's inequality for the class of quasi-subordinations.
\end{abstract}

\maketitle
\pagestyle{myheadings}
\markboth{S. Ponnusamy, R. Vijayakumar and K.-J. Wirths}{Bohr's Inequality}

\section{Introduction and Preliminaries}

Throughout the discussion, we let $\mathbb{D}=\{z\in \mathbb{C}:\,|z|<1\}$ denote the open unit disk
and $H^\infty $ denote the Banach algebra of all bounded analytic functions $f$ on the unit disk $\ID$
with the supremum norm $\|f\|_\infty :=\sup_{z\in \ID}|f(z)|$. Also, let
$${\mathcal B} = \{f\in H^\infty :\, \|f\|_\infty \leq 1 \}.
$$
By the maximum principle, the only members of ${\mathcal B}$ that touch the boundary $\partial\mathbb{D}$ of the unit disk are unimodular constant
functions. Thus, it is sometimes convenient to exclude constant functions (eg. in the discussion of subordination), but this does not affect our
discussion.

In 1914, Bohr \cite{Bohr-14} observed that if $f(z)=\sum_{n=0}^{\infty} a_n z^n  $ then
\be\label{PVW-eq0}
\sum_{n=0}^{\infty}|a_n|r^n\leq 1
\ee
holds for $0\leq r\leq 1/6$ and for all $f\in \mathcal{B}$. Then M. Riesz, I. Schur and F. W. Wiener independently showed that
\eqref{PVW-eq0} holds for $0\leq r\leq 1/3$ and that the constant $1/3$ is best possible. Other proofs and generalizations
 were also known in the literature, e.g.  Rizzonelli \cite{Rizz62}, Ricci \cite{RC1955}, Sidon \cite{sid}, Tomic \cite{tom},
and Paulsen et al  \cite{PPS}. In 1997, Boas and Khavinson  \cite{BoasKhavin-97-4} showed that a similar phenomenon occurs
for polydiscs in $\IC^n$. Other multidimensional variants of  Bohr-type theorems for bounded complete Reinhardt domains were obtained by
Aizenberg \cite{Aizen-00-1}. It is important to mention that Bombieri \cite{Bomb-1962} proved that
$$\sum\limits_{n=0}^{\infty} |a_n|\, r^n \leq \frac{3-\sqrt{8(1-r^2)}}{r} ~\mbox{ for }~ 1/3\leq r\leq 1/\sqrt{2}
$$
and for an alternate proof of this inequality, we refer to the recent article \cite{KayPon-AAA18}. So,  it is natural to
ask for the best constant $C(r)\geq 1$
such that
$$\sum\limits_{n=0}^{\infty} |a_n|\, r^n\leq C(r)
$$
whenever$f\in {\mathcal B}$. Later in \cite{BombBour-2004}, Bombieri and Bourgain proved that
$$\sum\limits_{n=0}^{\infty} |a_n|\, r^n< \frac{1}{\sqrt{1-r^{2}}} ~\mbox{ for $r> 1/\sqrt{2}$}
$$
which in turn shows that $C(r) \asymp (1-r^{2})^{-1/2}$ as $r\to 1$. In the same paper they also obtained a lower bound.

Besides these, several authors investigated the Bohr phenomenon in the recent years. For instance, Kayumov and Ponnusamy \cite{KayPon1,KP2017}
determined the Bohr radius for the class of
$p$-symmetric analytic  functions with multiple zeros at the origin, and introduced the notion of $p$-Bohr radius for harmonic functions and
obtained the $p$-Bohr radius for the class of odd harmonic bounded functions (see also \cite{ABS,KPS}) while in \cite{KayPon-AAA18}
the same authors discussed powered Bohr radius, originally discussed by Djakov and Ramanujan \cite{DjaRaman-2000}. Aytuna and Djakov  \cite{AD2013}
studied the Bohr property of bases for holomorphic functions, and
Ali et al. \cite{AAN2016}  discussed the Bohr radius for the class of starlike logharmonic mappings.
For further studies on the Bohr phenomenon, we refer to the survey articles \cite{AAP,GarMasRoss-2018} and
the references therein. On the other hand,  the authors in \cite{BaluCQ-2006} (see also Queff\'{e}lec \cite{Queff-1995}) extended
the work of Bohr to the setting of Dirichlet series.

In order to make the statements of the recent generalization and our present refined formulation, we need to introduce some basic notations.
Let ${\mathcal F}$ denote the set of all sequences  $\varphi =\{\varphi_{n}(r)\}_{n=0}^{\infty}$ of nonnegative continuous functions in $[0, 1)$
such that the series $\sum_{n=0}^{\infty} \varphi_{n}(r)$ converges locally uniformly with respect to $r \in [0, 1).$
For convenience, we let  $\Phi_N(r)=\sum_{n=N}^{\infty} \varphi_{n}(r)$ whenever $\varphi \in {\mathcal F}$.
Also, for $f(z)=\sum_{n=0}^{\infty} a_n z^n \in \mathcal{B}$ and $f_{0}(z):=f(z)-f(0),$ in what follows we let
$$ 
B_N(f, \varphi, r):= \sum_{n=N}^{\infty}|a_n| \varphi_{n}(r) \ \ \mbox{for}\ N \geq 0
$$
so that $B_0(f, \varphi, r)=|a_0|\varphi_{0}(r)+B_1(f, \varphi, r).$ In addition, we also let
$$A(f_0,\varphi, r):=\sum_{n=1}^{\infty}|a_n|^2 \left ( \frac{\varphi_{2n}(r)}{1+|a_{0}|}+\Phi_{2n+1}(r)\right)  ~\mbox{ and }~  \|f_{0}\|^2_{r}=\sum_{n=1}^{\infty}|a_n|^2 r^{2n}.
$$
In particular, when $\varphi_n(r)=r^n$, the formula for $A(f_0,\varphi, r)$ takes the following simple form
\beqq
A(f_0, r):=  \left(\frac{1}{1+|a_{0}|}+\frac{r}{1-r}\right)\|f_{0}\|^2_{r} ,
\eeqq
since $\Phi_{2n+1}(r)=r^{2n+1}/(1-r)$.


The following generalization is obtained recently by Kayumov et al. \cite{KaykhaPo}.

\begin{Thm}\label{PVW3-ThmB}
{\rm (\cite{KaykhaPo})}
Let $f\in {\mathcal B},$ $f(z)=\sum_{n=0}^{\infty} a_n z^n$ and $p \in (0,2].$ If $\varphi =\{\varphi_{n}(r)\}_{n=0}^{\infty}\in {\mathcal F}$ satisfies $p\varphi_{0}(r) > 2\Phi_1(r)$
for $r\in [0, R),$ where $R$ is the minimal positive root of the equation
$p\varphi_{0}(x)  = 2\Phi_1(x),$
then  the following sharp inequality holds:
\begin{equation}\label{KKP16-eq2}
	B_{f}(\varphi, p, r):= |a_{0}|^p \varphi_{0}(r)+ \sum_{n=1}^{\infty} |a_{n}| \varphi_{n}(r) \leq \varphi_{0}(r) ~\mbox{ for all $r \leq R$}.
\end{equation}
In the case when
$  \varphi_0(x) < (2/p)\sum_{k=1}^\infty \varphi_k(x)
$
in some interval $(R,R+\varepsilon)$, the number $R$ cannot be improved. If the functions $\varphi_k(x)$ ($k\geq 0$) are smooth functions then the last condition
is equivalent to the inequality
$$
\varphi_0'(R) <\frac{2}{p} \sum_{k=1}^\infty \varphi_k'(R)
$$
\end{Thm}

For convenience, we let $B_{f}(\varphi, r):=B_{f}(\varphi, 1, r)$, and it is natural to call this new majorant series as a {\em generalized majorant series} for
$f \in \mathcal{B}.$ We now begin with the discussion with the two basic properties of the new majorant series, which will be used in Section \ref{PVW3-sec6}.
This lemma is well-known in the fundamental case of $\varphi_{n}(r)=r^n$ (cf. \cite{PauSi})

\begin{lem}\label{KPV-lem1}
Let $f, g \in \mathcal{B}$ and $\varphi=\{\varphi_{n}(r)\}_{n=0}^{\infty}\in {\mathcal F}$. 
Then
\begin{itemize}
\item[(1)] $B_{f+g}(\varphi, r) \leq B_{f}(\varphi, r)+B_{g}(\varphi, r)$ for $r \in [0,1).$
\item[(2)]$B_{fg}(\varphi, r) \leq B_{f}(\varphi, r)B_{g}(\varphi, r)$ for $r \in [0,1)$ provided 
$\varphi_{k}$'s satisfy the additional condition
$\varphi_{m+n}(r)\leq \varphi_{m}(r)\varphi_{n}(r)$ for all $m, n \geq 0$ and $r \in [0,1).$
\end{itemize}	
Also, we note the trivial fact $B_{\alpha f}(\varphi, r)=|\alpha|B_{f}(\varphi, r)$ for all $\alpha \in \mathbb{C}.$
\end{lem}	
\begin{pf}
The proof of this lemma is easy, but for the sake of completeness, we include the proofs. Clearly,
$$ \sum_{n=0}^{\infty} |a_{n}+b_{n}| \varphi_{n}(r) \leq  \sum_{n=0}^{\infty} |a_{n}| \varphi_{n}(r)+\sum_{n=0}^{\infty} |b_{n}| \varphi_{n}(r).
$$
The product $fg$ takes the form
$$\sum_{n=0}^\infty c_n z^n =\sum_{n=0}^\infty a_n z^n \sum_{n=0}^\infty b_n z^n= \sum_{n=0}^\infty \left(\sum_{m+j=n}a_m b_j\right) z^n,
$$
which by equating the coefficients of $z^n$ on both sides gives
\beqq
c_n= \sum_{m+j=n}a_m b_j ~\mbox{ for each $n \ge 0$}.
\eeqq
Applying the triangle inequality to the last relation shows that
\beqq
\sum_{n=0}^\infty |c_n|\varphi_{n}(r)
&\le & \sum_{n=0}^\infty \left (\sum_{m+j=n}|a_m| \,|b_j| \right ) \varphi_{n}(r)\\
&\leq& \sum_{n=0}^\infty \sum_{m+j=n}|a_m|\varphi_{m}(r) |b_j|\varphi_{j}(r) \\
&= & \left (\sum_{m=0}^\infty |a_m|\varphi_{m}(r) \right )\sum_{j=0}^\infty|b_j|\varphi_{k}(r),
\eeqq
where in the second inequality above  we have used the inequality $\varphi_{m+j}(r)\leq \varphi_{m}(r)\varphi_{j}(r).$
The proof is complete.
\end{pf}

The article is organized as follows. In Section \ref{PVW3-sec2}, our aim is to improve Theorem~A 
(see Theorem \ref{PVW3-th1}) and
as a consequence, we establish few corollaries, remarks and examples. In Section \ref{PVW3-sec3}, we present an application which will provide
more examples of Bohr-type inequality in a refined form. Section \ref{PVW3-sec4} is dedicated to derive weighted Bohr radius for quasi-subordination
family of analytic functions whereas in Section \ref{PVW3-sec5} we deal with weighted Bohr type inequality for locally univalent  quasiconformal harmonic
mappings. Finally in Section \ref{PVW3-sec6}, using  the discussion of earlier sections, we consider the weighted Bohr type inequality for the derivative of Schwarz
functions along with few other related results.

\section{Weighted Bohr type inequality for functions in ${\mathcal B}$}\label{PVW3-sec2}

For the proof of our first theorem, we need the following lemma due to Carlson \cite{Carlson-40}.

\begin{Lem}\label{PVW1-lem2}
Suppose that $f\in {\mathcal B}$ and $f(z)=\sum_{n=0}^{\infty} a_n z^n$.
Then the following inequalities hold.
\begin{enumerate}
\item[{\rm (a)}] $|a_{2n+1}|\leq 1-|a_0|^2-\cdots - |a_n|^2,\ n=0,1,\ldots$ 
\item[{\rm (b)}] $|a_{2n}|\leq 1-|a_0|^2-\cdots -|a_{n-1}|^2 - \frac{|a_n|^2}{1+|a_0|} ,\ n=1,2,\ldots$. 
\end{enumerate}
Further, to have equality in {\rm (a)} it is necessary that $f$ is a rational function of the form
$$ f(z)=\frac{a_{0}+a_{1}z+ \cdots + a_{n}z^{n}+ \epsilon z^{2n+1}}{1+(\overline{a_n}z^n+ \cdots +\overline{a_0}z^{2n+1})\epsilon},~|\epsilon|=1,
$$
and to have equality in {\rm (b)} it is necessary that $f$ is a rational function of the form
$$ f(z)=\frac{a_{0}+a_{1}z+ \cdots + \frac{a_{n}}{1+|a_0|} z^{n}+ \epsilon z^{2n} }{1+\left(\frac{\overline{a_n}}{1+|a_0|}z^n+ \cdots +\overline{a_0}z^{2n}\right) \epsilon},
~|\epsilon|=1,
$$
where the term $a_0 \overline{a_n}^2 \epsilon$ is nonnegative real.
\end{Lem}


\begin{thm}\label{PVW3-th1}
Let $f\in {\mathcal B},$ $f(z)=\sum_{n=0}^{\infty} a_n z^n$ and $p \in (0,2].$ If $ \{\varphi_{n}(r)\}_{n=0}^{\infty} \in {\mathcal F}$ satisfies
the inequality
\be\label{PVW3-eq1}
p\varphi_{0}(r)  > 2\Phi_1(r),
\ee
then the following sharp inequality holds:
\be\label{PVW3-eq2}
 |a_{0}|^p \varphi_{0}(r)+ B_1(f, \varphi, r)+ A(f_0,\varphi, r)
\leq \varphi_{0}(r)~\mbox{for}~ r \leq R,
\ee
where $R$ is the minimal positive root of the equation
$	p\varphi_{0}(x)  = 2\Phi_1(x).	
$
In the case when $p\varphi_{0}(x)  < 2\Phi_1(x)$ in some interval $(R, R+\epsilon),$ the number $R$ cannot be improved.
\end{thm}
\begin{pf}
Using Carlson's lemma, we may rewrite $B_1(f, \varphi, r)$ as
\beq\label{PVW3-eq3}
B_1(f, \varphi, r)&= &\sum_{n=1}^{\infty}|a_{2n}|\varphi_{2n}(r)+\sum_{n=0}^{\infty}|a_{2n+1}|\varphi_{2n+1}(r)\nonumber\\
&\leq&  \sum_{n=1}^{\infty}\left[1-\sum_{k=0}^{n-1}|a_k|^2-\frac{|a_n|^2}{1+|a_0|}\right]\varphi_{2n}(r)
 	+\sum_{n=0}^{\infty}\left[1-\sum_{k=0}^{n}|a_k|^2\right]\varphi_{2n+1}(r)\nonumber\\
&=&   (1-|a_0|^2)\Phi_1(r)- \sum_{n=1}^{\infty} \frac{|a_n|^2}{1+|a_0|} \varphi_{2n}(r)- |a_1|^2 \Phi_3(r)-|a_2|^2 \Phi_5(r)-\cdots    \nonumber\\
&=&  (1-|a_0|^2)\Phi_1(r)- \sum_{n=1}^{\infty}|a_n|^2 \left[ \frac{\varphi_{2n}(r)}{1+|a_0|}+ \sum_{m=2n+1}^{\infty} \varphi_{m}(r)\right] \nonumber\\
&=& (1-|a_0|^2)\Phi_1(r)- A(f_0,\varphi, r)  \nonumber
\eeq
so that
\beqq
 |a_{0}|^p \varphi_{0}(r)+ B_1(f, \varphi, r)+ A(f_0,\varphi, r) \leq \Phi (p,|a_0|,r),
\eeqq
where
$$\Phi (p,|a_0|,r) = |a_0|^p \varphi_{0}(r)+(1-|a_0|^2)\Phi_1(r).
$$
It can be easily shown that for $0<p\leq 2$, the inequality $\Phi (p,|a_0|,r) \leq \varphi_{0}(r)$ holds for $r \leq R$,
under the condition \eqref{PVW3-eq1} (see also \cite[Theorem 1]{KaykhaPo}). This completes the proof of the inequality \eqref{PVW3-eq2}.

Now let us prove that $R$ is an optimal number. We consider the function
$$ f(z) =\frac{a-z}{1-a z}=a-(1-a^2)\sum\limits_{n=1}^\infty a^{n-1}z^n,\quad a\in [0,1).
$$
For this function, with $a_0=a$ and $a_n=-(1-a^2)a^{n-1}$, straightforward calculation shows that

\vspace{8pt}
$\ds  |a_{0}|^p \varphi_{0}(r)+ B_1(f, \varphi, r)+ A(f_0,\varphi, r)
$
\beqq
 &=&  a^p \varphi_{0}(r)+(1-a^2) \sum_{n=1}^{\infty} a^{n-1} \varphi_{n}(r)
 + (1-a^2)^2\sum_{n=1}^{\infty} a^{2n-2} \left[ \frac{\varphi_{2n}(r)}{1+a}+\Phi_{2n+1}(r)\right] \\
& =& \varphi_{0}(r)+(1-a)\left[2 \sum_{n=1}^{\infty}a^{n-1}\varphi_{n}(r)-p\varphi_{0}(r)\right]\\
&& \quad -(1-a)\left[(1-a) \sum_{n=1}^{\infty}a^{n-1}\varphi_{n}(r)+\left( \frac{1-a^p}{1-a}-p\right)\varphi_{0}(r)\right]\\
&&\qquad  + (1-a^2)^2\sum_{n=1}^{\infty} a^{2n-2} \left[ \frac{\varphi_{2n}(r)}{1+a}+\Phi_{2n+1}(r)\right] \\
&& =\varphi_{0}(r)+(1-a)\left[  2 \sum_{n=1}^{\infty}a^{n-1}\varphi_{n}(r)-p\varphi_{0}(r)\right]+ O((1-a)^2)
\eeqq
as $a\rightarrow 1^{-}.$ Now it is easy to see that the right hand side number in the above expression is $>\varphi_{0}(r)$ when $a$ is close to
$1$. The proof of the theorem is complete.
\end{pf}

In what follows, $ \floor{x} $ denotes the largest integer no more than $x,$ where $x$ is a real number.

\begin{cor}\label{PVW3-cor1}
Let $f \in \mathcal{B},$ $f(z)=\sum_{k=0}^{\infty}a_{k}z^{k}$, and $p \in (0,2].$
Then for each $n\in \IN$ with $n \geq 3,$ the following inequality holds:
\beqq
|a_{0}|^{p}+\sum_{k=1}^{n}|a_{k}|r^{k}+ \frac{1}{1+|a_0|} \sum_{k=1}^{s}|a_k|^2 r^{2k}+ \sum_{k=1}^{t}|a_k|^2 r^{2k+1}\left( \frac{1-r^{n-2k}}{1-r} \right) \leq 1
\eeqq
for $r \leq R_{n}(p),$ where $s=\floor{n/2},$ $t=\floor{(n-1)/2}$ and $R_{n}(p)$ is the minimal positive root of the equation
$$p(1-r)-2r(1-r^n)=0, ~\mbox{ or }~ p=2\sum_{k=1}^{n}r^k.
$$
\end{cor}
\begin{pf}
Set $\varphi_{k}(r)=r^k$ for $0\leq k \leq n,$ and $\varphi_{k}(r)=0$ for $k>n.$ First we remark that for the case $n=1,$ \eqref{PVW3-eq2} reduces  to
\beqq
|a_0|^p+ |a_1|r \leq 1~\mbox{for}~ r \leq R_1(p)=\frac{p}{2},
\eeqq
and, for the case $n=2,$ \eqref{PVW3-eq2} becomes
$$ |a_0|^p+|a_1|r+|a_2|r^2+|a_{1}|^2 \frac{r^2}{1+|a_0|} \leq 1~\mbox{for}~ r \leq R_2(p)=\frac{-1+\sqrt{1+2p}}{2},
$$
where $R_2(p)$ is the minimal positive root of the equation $ p=2r(1+r).$
Next we let $n \geq 3.$ Then \eqref{PVW3-eq2} is equivalent to
\beq\label{PVW3-eq6}
|a_0|^p+\sum_{k=1}^{n}|a_k|r^k+ I\leq 1~\mbox{ for all }~ r \leq R_n(p),
\eeq
 where  $R_{n}(p)$ is the minimal positive root of the equation $ p=2\sum_{k=1}^{n}r^k$ and
 $$
 I:=\sum_{k=1}^{s}|a_k|^2 \left[ \frac{r^{2k}}{1+|a_0|}+ r^{2k+1}+\cdots +r^{n} \right].
 $$
Now, the proof is divided into two cases. For the even values of $n \in \mathbb{N},$ we set $n=2m\ (m \geq 2).$ It follows easily that
\beqq
I &=& |a_1|^2 \left(\frac{r^{2}}{1+|a_0|}+r^3+\cdots +r^{2m}\right)\\
&& \quad +\cdots +|a_{m-1}|^2 \left(\frac{r^{2(m-1)}}{1+|a_0|}+r^{2m-1}+r^{2m}\right)+ |a_{m}|^2 \frac{r^{2m}}{1+|a_0|} \\
&=&  \frac{1}{1+|a_0|}\sum_{k=1}^{m}|a_k|^2 r^{2k} + \sum_{k=1}^{m-1}|a_k|^2 \left( r^{2k+1}+\cdots +r^{2m} \right)\\
&=& \frac{1}{1+|a_0|}\sum_{k=1}^{m}|a_k|^2 r^{2k} + \sum_{k=1}^{m-1}|a_k|^2 r^{2k+1}\left( \frac{1-r^{n-2k}}{1-r} \right)\\
\eeqq
so that \eqref{PVW3-eq6} gives the desired inequality  when $n\geq 4$ is even . For the odd values of $n \in \mathbb{N},$ we set $n=2m+1\ (m \geq 1).$ It follows that
\beqq
I &=& |a_1|^2 \left(\frac{r^{2}}{1+|a_0|}+r^3+\cdots +r^{2m+1}\right)\\
&& \quad +\cdots +|a_{m-1}|^2 \left(\frac{ r^{2(m-1)}}{1+|a_0|}+r^{2m-1}+r^{2m}+r^{2m+1}\right)+ |a_{m}|^2 \left(\frac{r^{2m} }{1+|a_0|}+r^{2m+1} \right)\\
&=& \frac{1}{1+|a_0|}\sum_{k=1}^{m}|a_k|^2 r^{2k} + \sum_{k=1}^{m}|a_k|^2 \left( r^{2k+1}+\cdots +r^{2m+1} \right)\\
&=& \frac{1}{1+|a_0|}\sum_{k=1}^{m}|a_k|^2 r^{2k} + \sum_{k=1}^{m}|a_k|^2 r^{2k+1}\left( \frac{1-r^{n-2k}}{1-r} \right).\\
\eeqq
Again, \eqref{PVW3-eq6} gives the desired inequality when $n\geq 3$ is odd. Combining the last two cases concludes the proof.
\end{pf}

Allowing $n\rightarrow \infty$, we obtain the following.

\beg\label{PVW3-eg1}
Suppose that $f \in \mathcal{B},$ $f(z)= \sum_{k=0}^{\infty} a_{k}z^{k}$ and $p\in (0,2].$ Then Theorem \ref{PVW3-th1} gives the following:
For $\varphi_{k}(r)=r^k \,(k \geq 0),$ we easily have
\begin{equation*}
|a_0|^p+ \sum_{k=1}^\infty |a_k|r^k+ \left(\frac{1}{1+|a_0|}+\frac{r}{1-r}\right) \sum_{k=1}^{\infty} |a_{k}|^2 r^{2k} \leq 1 ~\mbox{ for  }~ r \leq R(p)=\frac{p}{2+p}
\end{equation*}
and the constant $R(p)$ cannot be improved. The case $p=2$ is obtained in \cite[Theorem 2]{PoViWi-20}, the general case
is remarked in \cite[Remark 1]{PVW-2019},  and the inequality in this case does play a special role. The above inequality is a refined version of the
classical Bohr inequality.
\eeg


The values of $R_n(1)$ and $R_n(2)$ from Corollary \ref{PVW3-cor1} for certain choices of $n$ are listed in Table \ref{PVW3-tab1}.

\begin{table}[H]
\begin{center}
\begin{tabular}{|l|l||l|l||l|l||l|l||}
\hline
 {\bf $n$}   & $R_n(1)$ & $n$ &  $R_n(1)$  & $n$ & $R_n(2)$ & {\bf $n$}   & $R_n(2)$ \\
  \hline
2  & 0.366025  & 3 & 0.342508  & 2 & 0.618034   & 3 & 0.543689     \\
  \hline
4  & 0.336197  & 5 & 0. 334263   & 4 & 0.51879   & 5  & 0.50866     \\
  \hline
6  & 0.33364  & 7 & 0.333435  & 6 & 0.504138  & 7  & 0.502017     \\
  \hline
8  &0.333367   & 9 &0.333345   & 8 & 0.500994  & 9 & 0.500493    \\
 \hline
10  &0.333337   & 15 & 0.333333  & 10 & 0.500245  & 15 & 0.500008    \\
 \hline
20  & 0.333333  & 25 & 0.333333  & 20 & 0.5  & 25 & 0.5    \\
 \hline
30  & 0.333333   & 35 & 0.333333  & 30 & 0.5  & 35 & 0.5    \\
  \hline

\end{tabular}
\end{center}
\caption{$R_n(p)$ is the unique root of the equation
$p(1-r)- 2r(1-r^n) =0$  in $(0,1)$\label{PVW3-tab1}}
\end{table}

Note that when $p=1$ and $n\geq 16$, the approximate value of the root has no change in the $6^{th}$ decimal place and we
truncated the remaining digits. As $ n\rightarrow \infty$, the roots converge to $1/3$. Similarly, when $p=2$ and $n\geq 20$, the approximate
value of the root has no change in number up to  the $6^{th}$ decimal place and we truncated the remaining digits.
As $n \rightarrow \infty$, the roots converge to $1/2$. Moreover,
\beqq
R_2(1)&=& \frac{1}{2}(-1 + \sqrt{3}) \approx 0.366025 \\
R_3(1)&=&\frac{1}{3}\left(-1  - \frac{2 \times 2^{2/3}}{(41 + 3 \sqrt{201})^{1/3}} + \frac{(41 + 3 \sqrt{201})^{1/3}}{2^{2/3}} \right) \approx 0.342508 \\
R_2(2)&=& \frac{1}{2}(-1 + \sqrt{5}) \approx  0.618034 \\
R_3(2)&=& \frac{1}{3}\left(-1  - \frac{2 }{(17 + 3 \sqrt{33})^{1/3}} + (17 + 3 \sqrt{33})^{1/3} \right) \approx 0.543689.
\eeqq

\br
Moreover if $f \in \mathcal{B}$, then Corollary \ref{PVW3-cor1} for  $p=1$ shows that for each $n\geq 3$, the partial sum
$S_n(z,f)= \sum_{k=0}^{n} a_{k}z^{k}$ satisfies the sharp inequality
$$|S_n(z,f)| +\frac{1}{1+|a_0|} \sum_{k=1}^{s}|a_k|^2 r^{2k}+ \sum_{k=1}^{t}|a_k|^2 r^{2k+1}\left( \frac{1-r^{n-2k}}{1-r} \right)
\leq 1 ~\mbox{for}~ |z|=r \leq  R_{n}(1),
$$
where $s=\floor{n/2},$ $t=\floor{(n-1)/2}$ and $R_{n}(1)$ is the positive root of the equation 	
$$(1-r)-2r(1-r^n)=0, ~\mbox{ or }~ \sum_{k=1}^{n}r^k=\frac{1}{2}.
$$	
Recall from the previous corollary that $R_{1}(1)=\frac{1}{2}$ and $R_{\infty}(1)=\frac{1}{3}.$ Indeed, as we have to deal
the cases $n=1$ and $n=2$ separately, we have
\beqq
|S_1(z,f)|\leq |a_0|+ |a_1|r \leq 1~\mbox{for}~ r \leq R_1(1)=\frac{1}{2},
\eeqq
and
$$ |S_2(z,f)| + |a_{1}|^2 \frac{r^2}{1+|a_0|}\leq |a_0|+|a_1|r+|a_2|r^2+|a_{1}|^2 \frac{r^2}{1+|a_0|} \leq 1~\mbox{for}~ r \leq R_2(1)=\frac{-1+\sqrt{3}}{2}.
$$
We see that $R_{n}(1)$ is a decreasing function of $n$ from $\frac{1}{2}$
to $\frac{1}{3}.$ At this place it is worth pointing out that if $f \in \mathcal{B}$ then, according to Rogosinski \cite{Rogo-23}
(see also \cite{LanGair-86, SchurSzego-25}),
$$|S_n(z,f)|\leq 1 ~\mbox{for}~|z|\leq \frac{1}{2}.
$$
\er	

\section{An application  of Theorem \ref{PVW3-th1}} \label{PVW3-sec3}
Here is a simple application of Theorem \ref{PVW3-th1}. One can apply Theorem \ref{PVW3-th1} with
$\{\varphi_{n}(r)\}_{n=0}^{\infty}\in {\mathcal F}$, where $\varphi_{n}(r)=b_nr^n$ and $b_n\geq 0$ for $n \geq 0$.
For instance, the fractional derivative of $f$ of order $\alpha \in \mathbb{R}$ is defined by
$$ D^{\alpha} f(z)=\sum_{n=0}^{\infty} (n+1)^{\alpha}a_{n}z^{n}= (f*g)(z), \quad z\in \mathbb{D},
$$
where $g(z)=\sum_{n=0}^{\infty}(n+1)^{\alpha} z^n.$ Note that $D^{1}f(z)=(zf)'(z)$ and $D^{-1} f(z)=\frac{1}{z}\int_{0}^{z}f(t)dt.$
Then $(f_{0}*g_{0})(z)=(f*g)(z)-a_{0}$
and thus, by applying  Theorem \ref{PVW3-th1} with $\varphi_{0}(r)=1$ and $\varphi_{n}(r)=(n+1)^{\alpha} r^{n}$ for $n \geq 1$, we obtain
the following result.

\begin{thm}\label{PVW3-th2}
Suppose that $f(z)=\sum_{n=0}^{\infty}a_{n}z^{n}$ belongs to $\mathcal{B},$ $0<p \leq 2,$ and $\Phi_N(r) =\sum_{n=N}^{\infty}(n+1)^{\alpha}r^{n} $
with $\alpha \in \mathbb{R}$. Then
\begin{equation*}
|a_0|^p+ \sum_{n=1}^\infty (n+1)^\alpha |a_n|r^n+ \sum_{n=1}^{\infty}|a_n|^2 \left[ \frac{(2n+1)^{\alpha} r^{2n}}{1+|a_0|}+ \Phi_{2n+1}(r)\right] \leq 1 ~\mbox{ for  }~ r \leq R,
\end{equation*}
where $R$ is the minimal positive root of the equation $\Phi_1(r)={p}/{2}.$
\end{thm}


\begin{example}
If $f(z)=\sum_{n=0}^{\infty}a_{n}z^{n} \in \mathcal{B}$ and $0<p \leq2$, then
\begin{multline*}	
(1)~~|a_0|^p+ \sum_{n=1}^\infty (n+1) |a_n|r^n+ \left(\frac{1}{1+|a_0|}+\frac{r}{1-r}\right) \sum_{n=1}^{\infty}(2n+1) |a_{n}|^2 r^{2n} \\ +\frac{r}{(1-r)^2} \sum_{n=1}^{\infty}|a_{n}|^2 r^{2n}\leq 1~\mbox{ for  }~ r \leq R_{1}(p)=1-\sqrt{\frac{2}{2+p}}.
\end{multline*}
Indeed this follows from Theorem \ref{PVW3-th2} by setting $\varphi_{0}(r)=1$, $\varphi_{n}(r)=(n+1) r^{n}$ for $n \geq 1$
and  noting that
\beqq
\Phi_{2n+1}(r)  &=&  \sum_{k=2n+1}^{\infty} (k+1)r^{k}=  r^{2n+1}  \sum_{m=0}^{\infty} (m+1 +(2n+1))r^{m}\\
&=&r^{2n+1} \left ( \frac{1}{(1-r)^2}+\frac{2n+1}{1-r} \right ).
\eeqq
\begin{multline*}
 (2)~|a_0|^p+ \sum_{n=1}^\infty \frac{|a_n|}{n+1}r^n+ \sum_{n=1}^{\infty} |a_{n}|^2 \left[\frac{r^{2n}}{(2n+1)(1+|a_0|)}+ \frac{1}{r} \int_{0}^{r} \frac{t^{2n+1}}{1-t}\, dt \right] \leq 1
~\mbox{ for  }~ r \leq R_{2}(p),
\end{multline*}
where $R_{2}(p)$ is the unique positive root of the equation
$$-\frac{\log(1-r)}{r}=\frac{p+2}{2}, ~r \in (0,1).
$$
Indeed, by setting $\alpha =-1$, $\varphi_{0}(r)=1$ and $\varphi_{n}(r)=r^n/(n+1)$ for $n \geq 1$, we have at first
$$\Phi_1(r)= \sum_{n=1}^{\infty} \frac{r^{n}}{n+1}= -\frac{\log(1-r)}{r}-1,
$$
which is an increasing function of $r\in [0,1)$, and increases from	$0$ to $\infty.$ Secondly, it follows that
\beqq
\Phi_{2n+1}(r)  &=&  \sum_{k=2n+1}^{\infty} \frac{r^{k}}{k+1} = \sum_{m=1}^{\infty} \frac{r^{m+2n}}{m+2n+1}
=\frac{1}{r} \int_{0}^{r} \frac{t^{2n+1}}{1-t}  dt.
\eeqq
The desired conclusion follows from Theorem \ref{PVW3-th2}. By a standard computation (eg. using Mathematica),
we find that $ R_{2}(1)\approx 0.582812$ and $ R_{2}(2)\approx 0.796812$.
\begin{multline*}
(3)~|a_0|^p+ \sum_{n=1}^\infty (n+1)^2 |a_n| r^n+ \left(\frac{1}{1+|a_0|}+\frac{r}{1-r}\right) \sum_{n=1}^{\infty}(2n+1)^2 |a_{n}|^2 r^{2n} \\
+\sum_{n=1}^{\infty}|a_n|^2 r^{2n+1} \left[ \frac{1+r}{(1-r)^3}+ \frac{2(2n+1)}{(1-r)^2} \right]\leq 1
 \mbox{ for}~ r \leq R_{3}(p)=\frac{p+3-\sqrt{4p+9}}{p+2}.
\end{multline*}
Indeed, by setting $\alpha =2$, $\varphi_{0}(r)=1$ and $\varphi_{n}(r)=(n+1)^2 r^{n}$ for $n \geq 1$, we have
$$\Phi_1(r)= \sum_{n=1}^{\infty} (n+1)^2 r^{n} = \frac{1+r}{(1-r)^3}-1,
$$
and so the minimal positive root of the equation
$$\frac{1+r}{(1-r)^3}=\frac{p+2}{2}
$$
gives the number $R_{3}(p).$ Again, it follows from the choices of $\varphi_{n}(r)$'s that
\beqq
\Phi_{2n+1}(r)  &=&  \sum_{k=2n+1}^{\infty} (k+1)^2 r^{k}=  r^{2n+1}  \sum_{m=0}^{\infty} (m+1 +(2n+1))^2 r^{m}\\
&=&r^{2n+1} \sum_{m=0}^{\infty} \left[ (m+1)^2+2(m+1)(2n+1)+(2n+1)^2 \right] r^m \\
&=& r^{2n+1} \left[ \frac{1+r}{(1-r)^3}+\frac{2(2n+1)}{(1-r)^2}+\frac{(2n+1)^2}{1-r} \right].
\eeqq
The desired conclusion follows from Theorem \ref{PVW3-th2}. It is easy to find that $ R_{3}(1) = \frac{4-\sqrt{13}}{3} \approx 0.13148$ and $ R_{3}(2)= \frac{5-\sqrt{17}}{4} \approx 0.21922$.
\end{example}

\section{Weighted Bohr radius for quasi-subordination family}\label{PVW3-sec4}

Throughout this section, we let $\mathcal A$ denote the class of analytic functions in the unit disk $\ID$.
Unless otherwise stated, when we write $f,g\in \mathcal{A}$, we always assume the following power series representation :
\be\label{PVW3-eq1b}
f(z)=\sum_{k = 0}^\infty a_kz^k ~\mbox{ and }~g(z)=\sum_{k = 0}^\infty b_kz^k ~\mbox{ for $z\in \ID$}.
\ee

\bdefe {\rm \cite{Rob1970, Rob1970b}}
For any two analytic functions $f$ and $g$ in $\mathbb{D}$, we say that the function
$f$  is  quasi-subordinate to $g$  (relative to $W$),
denoted by $f \prec_{q} g$ in $\mathbb{D}$ if there exist two functions
$W\in\mathcal{B},\, \omega\in\mathcal{B}$ with $\omega(0)=0$ such that $f(z)=W(z) g(\omega(z))$.
\edefe

\begin{thm}\label{PVW3-qsth1}
Assume that $\varphi =\{\varphi_{n}(r)\}_{n=0}^{\infty}$ belongs to $\mathcal F$ such that $\varphi_{0}(r)=1$ and
\be\label{PVW-eq1a}
\varphi_{m+n}(r)\leq \varphi_{m}(r)\varphi_{n}(r) ~\mbox{ for all $m, n \geq 0$ and $r \in [0,1).$ }
\ee
If $f,g\in \mathcal{A}$ are given by \eqref{PVW3-eq1b} and $f\prec_q g$ in $\ID$, then we have
$$\sum_{k=0}^\infty |a_k|  \varphi_{k}(r) \le \sum_{k=0}^\infty |b_k| \varphi_{k}(r) \  \ \text{for all}\  \ r \le R,
$$
where $R$ is the minimal positive root of the equation  $1 = 2\Phi_1(x)$,
$\Phi_1(x) =\sum_{n=1}^{\infty} \varphi_{n}(x).
$
\end{thm}
\begin{pf}
We remark that this theorem was proved in \cite[Theorem 2.1]{AlkKayPon} for $\varphi_{k}(r)=r^{k}~ (k\geq 0).$
We follow the method of proof of \cite[Theorem 2.1]{AlkKayPon}. Suppose that $f\prec_q g$. Then there exist two analytic
functions $W$ and $\omega$ satisfying $\omega(0)=0$, $|\omega(z)|\le1$ and
$|W(z)|\le1$ for all $z\in \ID$ such that 
\begin{equation}\label{Eq3}
f(z)=W(z)g(\omega(z)).
\end{equation}
Now for the analytic function $\omega(z)=\sum_{n=1}^\infty \alpha_n z^n,$   the
$k$-th power of $\omega$, where $k \in \IN,$ can be written as
\begin{equation}\label{Eq4}
\omega^k(z) = \sum_{n=k}^\infty \alpha_n^{(k)} z^n 
=z^k \sum_{m=0}^\infty \alpha_{m+k}^{(k)}z^{m}.
\end{equation}
For $k=0$, we set $\omega^0(z)=1$ so that $\alpha_0^{(0)}=1$ and $\alpha_n^{(0)}= 0$ for $n\geq 1$. As $\omega \in {\mathcal B}$ with $\omega(0)=0$, we may
write $\omega^k(z)=z^k\omega_1(z)$ with  $\omega_1 \in {\mathcal B}$. Applying Theorem \ref{PVW3-th1} to $\omega_1$ with $p=1$, it follows in particular that
\be\label{EQ5}
\sum_{n=k}^\infty|\alpha_n^{(k)}|\varphi_{n-k}(r)\le \varphi_{0}(r)=1 \ \ \text{for all}\ \ r\le R .
\ee
For the analytic function $W(z)$, we may write $W(z)=\sum_{m=0}^\infty w_m z^m$ and thus, by Theorem \ref{PVW3-th1}, we have
\be\label{EQ6}
\sum_{m=0}^\infty |w_m|\varphi_{m}(r)\le1 \ \ \text{for all}\ \ r\le R .
\ee
The relation \eqref{Eq3} with the help of \eqref{Eq4} takes the form (cf. \cite[Theorem 2.1]{AlkKayPon})
$$\sum_{k=0}^\infty a_k z^k =\sum_{m=0}^\infty w_m z^m \sum_{k=0}^\infty B_k z^k= \sum_{k=0}^\infty \left(\sum_{m+j=k}w_m B_j\right) z^k,
$$
which by equating the coefficients of $z^k$ on both sides gives 
\begin{equation}\label{Eq7}
a_k = \sum_{m+j=k}w_m B_j ~\mbox{ for each $k \ge 0$},
\end{equation}
where $B_k = \sum_{n=0}^k b_n\alpha_k^{(n)}$.
Applying the triangle inequality to the last relation shows that
\beqq
\sum_{k=0}^\infty |a_k|\varphi_{k}(r)
&\le & \sum_{k=0}^\infty \left (\sum_{m+j=k}|w_m| \,|B_j| \right ) \varphi_{k}(r)\\
&\leq &\sum_{k=0}^\infty \sum_{m+j=k}|w_m|\varphi_{m}(r) |B_j|\varphi_{j}(r) \\
&= & \left (\sum_{m=0}^\infty |w_m|\varphi_{m}(r) \right )\sum_{k=0}^\infty|B_k|\varphi_{k}(r)\\
&\le& \sum_{k=0}^\infty|B_k|\varphi_{k}(r) ~\mbox{ for all $r\le R$ \quad (by \eqref{EQ6})},
\eeqq
where in the second inequality above  we have used the inequality \eqref{PVW-eq1a}.
Also, because $|B_k| \leq  \sum_{n=0}^k |b_n|\, |\alpha_k^{(n)}|$, we obtain that
\beqq
\sum_{k=0}^\infty|B_k|\varphi_{k}(r)
&\le & \sum_{k=0}^\infty\sum_{n=0}^k |b_n|\, |\alpha_k^{(n)}|\varphi_{k}(r)
=\sum_{k=0}^\infty|b_k| \sum_{n=k}^\infty|\alpha_n^{(k)}|\varphi_{n}(r) \\
&\leq& \sum_{k=0}^\infty|b_k| \left (\sum_{n=k}^\infty|\alpha_n^{(k)}|\varphi_{n-k}(r)\right ) \varphi_{k}(r)
\quad \mbox{(by \eqref{PVW-eq1a})}\\
&\le& \sum_{k=0}^\infty|b_k|\varphi_{k}(r)  ~\mbox{ for all $r\le R$ \quad (by \eqref{EQ5})}
\eeqq
and hence, we deduce that
$$\sum_{k=0}^\infty |a_k|\varphi_{k}(r) \le \sum_{k=0}^\infty|B_k|\varphi_{k}(r)
\le \sum_{k=0}^\infty|b_k|\varphi_{k}(r) \ \ \text{for all}\ \ r\le R.
$$
The proof of Theorem \ref{PVW3-qsth1} is complete.
\end{pf}
%
%

\br
Clearly, the conclusion of Theorem \ref{PVW3-qsth1} continues to hold if the assumption $f\prec_q g$ is replaced by either $f\prec g$ or
the majorization condition  $|f(z)|\leq |g(z)|$ in $\ID$.
\er

The following theorem is due to Rogosinski \cite{Rogo-43}.

\begin{Thm}[Rogosinski's Theorem]\label{PVW-TheRog1}
Suppose that $f,g\in \mathcal{A}$ are given by \eqref{PVW3-eq1b},  $A_n=\sum_{k=1}^n |a_k|^2$ and $B_n=\sum_{k=1}^n |b_k|^2 $.
If $f\prec g$ in $\ID,$ then $A_n \le B_n$ for each $ n\geq 1.$
\end{Thm}

Goluzin \cite{Goluzin} observed that Rogosinski's theorem (with $\psi_{k}(r)=r^k$) can be applied to obtain the following more general result,
which has some interesting consequences. Here we state a general result which can be applied to a variety of situations.

\bthm[Simple generalization of Goluzin's Lemma]
\label{PVW-ThGolu1}
Suppose that $f,g\in \mathcal{A}$ are given by \eqref{PVW3-eq1b}, $f\prec g$ in $\ID,$ and  $\{\psi_{n}(r)\}_{n=1}^{\infty}$ is a decreasing sequence of
nonnegative functions in $[0, r_\psi)$. Then
$$
\sum_{k=1}^\infty |a_k|^2 \psi_{k}(r) \leq \sum_{k=1}^\infty |b_k|^2 \psi_{k}(r)\ \ \text{for}\ \ r\in [0,r_\psi).
$$
\ethm
\begin{pf}
Following the method of proof of \cite[Theorem 6.3]{Duren} and the notation of Rogosinski's
Theorem~C 
(i.e.  $A_n \le B_n$ for each $n \geq1 $), a summation by parts gives the desired inequality. So, we omit the details.
\end{pf}


\begin{thm}\label{PVW3-qsth2}
Assume the hypotheses of Theorem \ref{PVW3-qsth1} and, in addition, suppose that $\{\psi_k(r)\}_{k\geq 1}$ is a decreasing sequence
of non-negative functions defined in $[0, r_\psi)$. Then we have
$$ \sum_{k=0}^\infty |a_k|  \varphi_{k}(r)+ \lambda(r) \sum_{k=1}^\infty |a_k|^2 \psi_k(r)  \le
\sum_{k=0}^\infty |b_k| \varphi_{k}(r) +\lambda(r) \sum_{k=1}^\infty |b_k|^2 \psi_k(r)
$$
holds for all  $r \le \min \{R,r_\psi\}$, where $R$ is as in Theorem~\ref{PVW3-qsth1} and $\lambda (r)$ is a non-negative function of $r$ defined in $[0,1]$.
%
\end{thm}
\bpf
According to Theorem \ref{PVW3-qsth1}, we obtain from the assumptions that
\be\label{Eq9}
\sum_{k=0}^\infty |a_k|  \varphi_{k}(r) \le \sum_{k=0}^\infty |b_k| \varphi_{k}(r) \  \ \text{for all}\  \ r \le R.
\ee
Finally, by \eqref{Eq3}, it is also clear that
$$|f(z)|^2 \leq |g(\omega(z))|^2 ~\mbox{ for }~z\in \ID
$$
and thus, as in the proof of Rogosinski's Theorem \cite{Rogo-43}, we can easily obtain that
(see Robertson \cite{Rob1970} and  Pommerenke \cite[Theorem 2.2]{P})
\be\label{Eq10}
\sum_{k=1}^n |a_k|^2 \le \sum_{k=1}^n |b_k|^2 \  \ \text{for}\  \ n=1,2,\ldots.
\ee
Applying Theorem \ref{PVW-ThGolu1},  we obtain that
\be\label{Eq11}
\sum_{k=1}^\infty |a_k|^2 \psi_k(r) \leq \sum_{k=1}^\infty |b_k|^2 \psi_k (r)\ \ \text{for}\ \ r\in [0,r_\psi).
\ee
The desired inequality follows from \eqref{Eq9} and \eqref{Eq11}.
\epf

\br\label{PVW3-qsrem1}
Theorem \ref{PVW3-qsth1} holds for $f\prec g$ in $\ID$ (instead of $f\prec_q g$ in $\ID$). In this case, Theorem \ref{PVW3-qsth2} with
$\varphi_{k}(r)= r^{k} \,(k\geq 0)$,  $\psi_{k}(r)= r^{2k}\,(k\geq 1)$ and $\lambda (r) = \frac{1}{1+|a_0|}+\Phi_1(r)$, is well-known from
\cite[Lemma 2]{PVW-2019}. In this choice the value of $\min \{R,r_\psi\}$ turns out to be $1/3$. It is worth remarking that
the weight function $\lambda (r) = \frac{1}{1+|a_0|}+\Phi_1(r)$ appears very naturally in some specific situation which refines
the classical Bohr inequality and is due to Carlson lemma (see Lemma~B). 
For details, revisit Example \ref{PVW3-eg1} and  \cite{PoViWi-20,PVW-2019}.
\er

%

\begin{thm} \label{PVW3-qsth3}
Under the hypothesis of Theorem \ref{PVW3-qsth2}, we have the following inequality
$$ \sum_{k=0}^\infty |a_k|  \varphi_{k}(r)+\lambda (r)\sum_{k=1}^\infty |a_k|^2 k r^{2k} \le
\sum_{k=0}^\infty |b_k| \varphi_{k}(r) + \lambda (r)\sum_{k=1}^\infty |b_k|^2 k r^{2k}
$$
holds for all  $r \le \min\{R,\frac{1}{\sqrt{2}}\}.$

%

In particular, for $f\prec g$ and $\varphi_{k}(r)= r^{k} \,(k\geq 0),$  the inequality
$$\sum_{k=0}^\infty |a_k|r^k+\lambda (r)\sum_{k=1}^\infty |a_k|^2 k r^{2k} \le
 \sum_{k=0}^\infty |b_k|r^k +\lambda (r)\sum_{k=1}^\infty |b_k|^2 k r^{2k}
$$
%
holds for all  $r\leq  1/3$.
\end{thm}
\begin{pf}
 We set $\psi_{k}=k r^{2k}$ for $k\in \mathbb{N}$.
Then the sequence $\{k r^{2k}\}$ is non-increasing if and only if
$$ r \leq \left(\frac{k}{k+1}\right)^{1/2}, \ \ k=1,2,\ldots ,
$$
which holds for all $k$ if it holds for $k=1$. This gives the condition  $r \leq r_\psi =\frac{1}{\sqrt{2}}.$
The desired inequality follows from the method of proof Theorem \ref{PVW3-qsth2}.
\end{pf}

\begin{cor}
Under the hypothesis of Theorem \ref{PVW3-qsth2}, we have the following inequality
$$\sum_{k=0}^\infty |a_k|  \varphi_{k}(r)+\lambda (r)\sum_{k=1}^\infty |a_k|^2 k^2 r^{2(k-1)} \le
\sum_{k=0}^\infty |b_k|\varphi_{k}(r)+ \lambda (r)\sum_{k=1}^\infty |b_k|^2 k^2 r^{2(k-1)}
$$
holds for all  $r \le \min\{R,\frac{1}{2}\},$  and $R$ is as in Theorem \ref{PVW3-qsth1}.

For $f\prec g$ and $\varphi_{k}(r)= r^{k}\, (k\geq 0),$  the  inequality	
$$\sum_{k=0}^\infty |a_k|r^k+\lambda (r)\sum_{k=1}^\infty |a_k|^2 k^2 r^{2(k-1)} \le
\sum_{k=0}^\infty |b_k|r^k +\lambda (r)\sum_{k=1}^\infty |b_k|^2 k^2 r^{2(k-1)}
$$
holds for all  $r\leq  1/3$.
\end{cor}
\begin{pf}
Following the method of proof Theorem \ref{PVW3-qsth2},   we set $\psi_{k}=k^2 r^{2(k-1)}$ for $k\in \mathbb{N}$.
It can be easily seen that $\psi_{k}$ is decreasing if and only if $r\leq k/(k+1)$ and thus,  the sequence $\{\psi_{k}\}$ is decreasing for all $k\in \IN$,
provided $r\leq 1/2$. The conclusion follows.
\end{pf}

We would like to emphasize that for a fixed $g$ in the assumption $f\prec_q g$, fixed $\lambda(r)$, $\varphi_{k}(r)$ and $\psi_k(r)$ $(k\geq 0)$, it is possible to
find an interval for $r\in (0,1)$ such that $\sum_{k=0}^\infty |b_k| \varphi_{k}(r) +\lambda(r) \sum_{k=1}^\infty |b_k|^2 \psi_k(r)\leq 1$. By doing so, the
largest value of $r$ satisfying the last inequality gives the Bohr radius with Bohr-type inequality in a more general setting. Our approach in the above theorems
and corollaries is to provide a method of obtaining more such Bohr-type inequalities. However the past known examples are obtained by fixing $\varphi_{k}(r)=r^k$.

\section{Bohr radius for locally univalent harmonic mappings}\label{PVW3-sec5}

A sense-preserving harmonic mappings $f$ of the form $f=h+\overline{g},$ where $h$ and $g$ are analytic in $\ID$, is said to be
$K$-quasiconformal if $|g'(z)|\leq k |h'(z)|$ in the unit disk, for $k=\frac{K-1}{K+1} \in [0,1]$.  See \cite{KPS} for discussion on
Bohr radius for quasiconformal harmonic mappings.

\begin{lem}\label{PVW3-lem3.1}
Let $\{\psi_{n}(r)\}_{n=1}^{\infty}$ be a decreasing sequence of nonnegative functions in $[0,r_\psi)$, and
$g, h\in \mathcal{A}$ such that $|g'(z)|\leq k |h'(z)|$ in $\ID$ and for some $k \in [0,1],$ where $h(z)=\sum_{n=0}^\infty a_nz^n$
and $g(z)=\sum_{n=0}^\infty b_nz^n$. Then
$$\sum_{n=1}^\infty |b_n|^2 \psi_{n}(r) \leq k^2 \sum_{n=1}^\infty |a_n|^2 \psi_{n}(r) \ \ \text{for}\ \ r\in [0,r_\psi).
$$
\end{lem}
\begin{pf}
Since $|g'(z)|\leq k |h'(z)|$ in $\ID$ and thus, as in the proof of Rogosinski's Theorem \cite{Rogo-43}, we can easily obtain that
\beqq
\sum_{m=1}^n m^2 |b_m|^2 r^{2(m-1)} \le k^2 \sum_{m=1}^n m^2 |a_m|^2 r^{2(m-1)} \  \ \text{for}\  \ n=1,2,\ldots.
\eeqq
We integrate the last inequality with respect to $r^2$ and obtain
\beqq
\sum_{m=1}^n m |b_m|^2 r^{2m} \le k^2 \sum_{m=1}^n m|a_m|^2 r^{2m} \  \ \text{for}\  \ n=1,2,\ldots.
\eeqq
One more integration (after dividing by $r^2$) gives
\beqq
\sum_{m=1}^n  |b_m|^2  r^{2m} \le k^2 \sum_{m=1}^n |a_m|^2 r^{2m}  \  \ \text{for}\  \ n=1,2,\ldots.
\eeqq
 By letting $r$ tends to $1,$ one has
\beqq
\sum_{m=1}^n  |b_m|^2  \le k^2 \sum_{m=1}^n |a_m|^2  \  \ \text{for}\  \ n=1,2,\ldots.
\eeqq
Now, applying Goluzin lemma (see Theorem \ref{PVW-ThGolu1}) for the given set $\{\psi_{n}(r)\}$ gives the desired result.
\end{pf}


\begin{thm}\label{PVW3-thm3.2}
Assume that $g, h\in \mathcal{A}$ such that $|g'(z)|\leq k |h'(z)|$ in $\ID$ and for some $k \in [0,1],$ where
$g(z)=\sum_{n=0}^\infty b_nz^n$, and $h(z)=\sum_{n=0}^\infty a_nz^n$ satisfies the condition ${\rm Re}\,h(z)\leq 1$ in
$\ID$ and $h(0)=a_{0}$ is positive. If $ \{\varphi_{n}(r)\}_{n=1}^{\infty} \in {\mathcal F}$ is a decreasing sequence,
where $\Phi_1(r)= \sum_{n=1}^{\infty} \varphi_{n}(r)$, and satisfies the inequality
\beq\label{Eq 4}
1 >\frac{2}{p}(1+k)\Phi_1(r),
\eeq
for some $p\in (0,1]$, then the following sharp inequality holds:
\beq\label{Eq 5}
a_{0}^p+\sum_{n=1}^\infty |a_n| \varphi_{n}(r)+\sum_{n=1}^\infty |b_n| \varphi_{n}(r) \leq 1 \  \ \text{for all}\  \ r \le R,
\eeq
where $0<p\leq 1$ and $R$ is the minimal positive root of the equation
$$ 1=\frac{2}{p}(1+k)\Phi_1 (x).
$$
In the case when $1<\frac{2}{p}(1+k)\Phi_1(x)$	in some interval $(R, R+\epsilon),$ the number $R$ cannot be improved.
\end{thm}	
\begin{pf}
We recall that if $P(z)=\sum_{n=0}^{\infty}p_{n}z^{n}$ is analytic in $\ID$ such that ${\rm Re}\,p(z)>0$ in $\ID,$ then
$|p_{n}|\leq 2 {\rm Re}\,p_{0}$ for all $n \geq 1.$ Applying this result to $p(z)=1-h(z)$ leads to
$|a_{n}|\leq 2(1-a_{0})$ for all $n \geq 1.$ Thus, as in the proof of Theorem \ref{PVW3-th1}, we can easily obtain from Lemma \ref{PVW3-lem3.1} that	
\beqq	
\sum_{n=1}^\infty |b_n|^2 \varphi_{n}(r) \leq k^2 \sum_{n=1}^\infty |a_n|^2 \varphi_{n}(r) \leq 4k^2 (1-a_{0})^2 \sum_{n=1}^{\infty} \varphi_{n}(r)
=4k^2 (1-a_{0})^2\Phi_1(r).
\eeqq	
Consequently, it follows from the classical Schwarz inequality that
\beqq	
\sum_{n=1}^\infty |b_n| \varphi_{n}(r) \leq  \sqrt{\sum_{n=1}^\infty |b_n|^2 \varphi_{n}(r)} \sqrt{\sum_{n=1}^{\infty}\varphi_{n}(r)}
\leq 2k (1-a_{0})\Phi_1(r)
\eeqq
so that	
\beq
a_{0}^p+\sum_{n=1}^\infty |a_n| \varphi_{n}(r)+\sum_{n=1}^\infty |b_n| \varphi_{n}(r)
&\leq& a_{0}^p+ 2(1-a_{0})(1+k) \Phi_1(r) \nonumber \\
&=& 1+(1-a_{0}) \left[2(1+k)\Phi_1(r)-\left(\frac{1-a_{0}^p}{1-a_{0}}\right) \right] \nonumber \\
&\leq& 1, \ \ \text{by Eqn. \eqref{Eq 4}} \nonumber,
\eeq	
for all $r \leq R,$ by the definition of $R.$ In the last step, we have used the fact that the function
$$B(x)=\frac{1-x^p}{1-x}, \quad x\in [0,1)
$$
is a decreasing function of $x\in [0,1)$ so that
$$ B(x) \geq \lim_{x\rightarrow 1^{-}}\frac{1-x^p}{1-x}\, =p.
$$
This proves the desired inequality \eqref{Eq 5}.
Moreover, sharpness can be seen by considering functions as in Theorem \ref{PVW3-th1}.
\end{pf}	
	
\begin{example}
Suppose that $f(z)= h(z)+\overline{g(z)}=\sum_{n=0}^\infty a_nz^n +\overline{\sum_{n=1}^\infty b_nz^n}$ is a sense-preserving $K$-quasiconformal harmonic mapping of the disk $\ID,$ where $h$ satisfies the condition ${\rm Re}\,h(z)\leq 1$ in $\ID$ and $h(0)=a_{0}$ is positive. Then, by choosing $\varphi_{k}(r)= r^{k} \,(k\geq 0)$ in Theorem \ref{PVW3-thm3.2} and $p=1$, we obtain the following sharp inequality (see \cite[Theorem 1.3]{KPS})
\beqq
a_{0}+\sum_{n=1}^\infty |a_n| r^n+\sum_{n=1}^\infty |b_n| r^n \leq 1 \  \ \text{for all}\  \ r \le \frac{K+1}{5K+1}.
\eeqq
\end{example}

\section{Weighted Bohr for the derivative of Schwarz functions}\label{PVW3-sec6}

The following three theorems extend the work of Bhowmik and Das \cite[Theorems 1, 2 and 3]{BowDas} in our general setting.

\begin{thm}\label{thm 1}
Let $f$ be a Schwarz function, i.e. $f\in {\mathcal B}$ such that $f(0)$. If $\varphi=\{\varphi_{n}(r)\}_{n=0}^{\infty}\in {\mathcal F}$ such that
\be\label{PoKV_eq1}
\varphi_0(r) \geq  2\sum_{n=1}^{\infty} (n+1)\varphi_{n}(r),
\ee
then the following sharp inequality holds:
\be\label{PoKV_eq2}
B_{f'}(\varphi,r)
\leq \varphi_0(r)~\mbox{for all}~ r \leq R_{0},
\ee
where $R_{0}$ is the minimal positive root of the equation
$$ \varphi_0(x) = 2\sum_{n=1}^{\infty} (n+1)\varphi_{n}(x).
$$
In the case when $\varphi_0(x)< 2\sum_{n=1}^{\infty} (n+1)\varphi_{n}(x)$ in some interval $(R_{0}, R_{0}+\epsilon),$ the number $R_{0}$ cannot be improved.
\end{thm}
\begin{pf}
Let $f(z)=\sum_{n=1}^{\infty} a_{n}z^n = z\sum_{n=0}^\infty a_{n+1}z^n.$ Then, $f'(z)=\sum_{n=0}^{\infty}(n+1)a_{n+1}z^n$ and
$$
B_{f'}(\varphi, r) =\sum_{n=0}^{\infty}(n+1)|a_{n+1}|\varphi_{n}(r).
$$
Using Wiener's estimates  $|a_{n+1}| \leq 1-|a_{1}|^2\leq 2(1-|a_{1}|)$ for $n \geq 1$, we obtain that
\beq
B_{f'}(\varphi, r) 
&\leq& \varphi_0(r)+(1-|a_{1}|) \left[2\sum_{n=1}^{\infty} (n+1)\varphi_{n}(r) -\varphi_0(r) \right] \nonumber \\
&\leq& \varphi_0(r), \ \ \text{by Eqn. \eqref{PoKV_eq1}} \nonumber,
\eeq	
for all $r \leq R_{0},$ by the definition of $R_{0}.$ This
completes the proof of the inequality \eqref{PoKV_eq2}. Now let us prove that $R_{0}$ is an optimal number. We consider the function $f=\varphi_a$ given by
$$
\varphi_a (z)  = z\left(\frac{a-z}{1-az}\right) =az - (1-a^2)\sum_{n=1}^\infty a^{n-1} z^{n+1}, \ z\in\ID,
$$
where $a\in (0,1)$. For this function, 
straightforward calculations show that
\beqq
B_{f'}(\varphi, r)
& =& \varphi_0(r)+(1-a)\left[  2 \sum_{n=1}^{\infty} a^{n-1} (n+1)\varphi_{n}(r) -\varphi_0(r)\right] \\	
&& \quad -(1-a)^2\left[\sum_{n=1}^{\infty} a^{n-1} (n+1)\varphi_{n}(r) \right]\\
& =& \varphi_0(r)+(1-a)\left[  2 \sum_{n=1}^{\infty} a^{n-1} (n+1)\varphi_{n}(r) -\varphi_0(r)\right]+ O((1-a)^2)
\eeqq
as $a\rightarrow 1^{-}.$  Now it is easy to see that the right hand side  is $>1$ when $a$ is close to $1$.
The proof of the theorem is complete.	
\end{pf}

\br
The choice $\varphi_{k}(r)=r^{k}~ (k\geq 0)$ in Theorem \ref{thm 1} gives \cite[Theorem 1]{BowDas}, where the corresponding value of $R_0$ is $1-\sqrt{2/3}.$ 	
\er

\begin{thm}\label{thm 2}
Assume that $\varphi_{k}$'s are defined as in Theorem \ref{PVW3-qsth1}. Also, let $f,g\in \mathcal{A}$. Then we have the following:
\begin{itemize}
\item [{\rm (a)}] if $f\prec g$ in $\ID$ then 	$B_{f'}(\varphi, r) \leq 	B_{g'}(\varphi, r)$ for $r \leq r_{0}=\min\{R,R_0\}.$
\item[{\rm (b)}] if $|f(z)| \leq |g(z)|$ for all $z \in \mathbb{D},$ and $g(0)=0$ then $B_{f'}(\varphi, r) \leq 	B_{g'}(\varphi, r)$ for $r \leq R_{0},$
\end{itemize}
where $R$ and  $R_{0}$ are as in Theorems \ref{PVW3-qsth1} and  \ref{thm 1}, respectively.
\end{thm}	
\begin{pf}
To prove the first part, we let $f\prec g$, Then $f(z)=g(w(z))$ for $z \in \mathbb{D}$, where $w$ is a Schwarz function.
Since $f'(z)=w'(z)g'(w(z)),$  it follows that	
$$B_{f'}(\varphi, r) \leq B_{w'}(\varphi, r)B_{g' \circ w}(\varphi, r).
$$
As $g' \circ w \prec g',$ the remark followed by Theorem \ref{PVW3-qsth1}, implies that
$$B_{g' \circ w}(\varphi, r) \leq B_{g'}(\varphi, r) ~\mbox{ for $r \leq R,$}
$$
and from Theorem \ref{thm 1}, we have  $B_{w'}(\varphi, r) \leq 1$ for $r \leq R_{0}.$
Combining the last two observations gives the desired inequality (a).

To prove the second inequality (b), it suffices to
assume next that $|f(z)|<|g(z)|$ for all $z \in \mathbb{D}\setminus\{0\}.$ Then there exists an analytic self map
$h$ of $\mathbb{D}$ such that $f(z)=h(z)g(z)$ for $z \in \mathbb{D}.$ Now $f'(z)=h'(z)g(z)+h(z)g'(z).$
Further, observing that $B_{g/z}(\varphi, r)\leq	B_{g'}(\varphi, r)$ provided $g(0)=0,$ we have
\beqq
B_{f'}(\varphi, r)
&\leq& \left(B_{zh'}(\varphi, r)+B_{h}(\varphi, r)\right)B_{g'}(\varphi, r).
\eeqq		
Following similar lines of calculations as in the proof of Theorem \ref{thm 1}, it can be shown that
$$ B_{zh'}(\varphi, r)+B_{h}(\varphi, r) \leq 1  ~\mbox{ for }~ r \leq R_{0},
$$
and hence $B_{f'}(\varphi, r) \leq 	B_{g'}(\varphi, r)$ for $r \leq R_{0}.$	
\end{pf}
	

\br
If we choose $\varphi_{k}(r)=r^{k}$ $(k\geq 0)$, then Theorem \ref{thm 2} gives \cite[Theorem 2]{BowDas}, where the corresponding value of $r_0$ is $1-\sqrt{2/3}.$
\er

\begin{thm}	\label{PVW2-the10}
Assume that $\varphi_{k}$'s are defined as in Theorem \ref{PVW3-qsth1}. Let $f(z)=\sum_{n=0}^{\infty}a_{2n+1}z^{2n+1}$ and
$g(z)=\sum_{n=0}^{\infty}b_{2n+1}z^{2n+1}$ be two analytic functions defined on $\mathbb{D}$ such that $|f(z)| \leq |g(z)|$
for all $z \in \mathbb{D}.$ If
$$1> 2\sum_{n=1}^{\infty} \varphi_{2n}(r),
$$
then the following inequality holds:
\be\label{PoKV_eq4}
B_{f}(\varphi, r)\leq B_{g}(\varphi, r) ~\mbox{for all}~ r \leq R,
\ee
where $R$ is the minimal positive root of the equation
$ 	1= 2\sum_{n=1}^{\infty} \varphi_{2n}(x).
$
\end{thm}
\begin{pf}
As the method of proof is based on the proofs of Theorem \ref{thm 2} and \cite[Theorem 3]{BowDas}, we omit the details.
\end{pf}	

\br
If we choose $\varphi_{k}(r)=r^{k}$ $(k\geq 0)$, then Theorem \ref{PVW2-the10} gives \cite[Theorem 3]{BowDas}, where the corresponding value of $R$ is $\sqrt{1/3}.$
\er

\end{document}